\documentclass{article}
\usepackage{amssymb}
\usepackage{amsmath}
\usepackage{epsf}
\usepackage{color}
\def \C {{\Bbb C\,}}
\def \CC {\overline{\Bbb C\,}}
\def\R {{\Bbb R}}

\def\boF{{\cal F}}
\def\boH{{\cal H}}

\def\boZ{{\cal Z}}
\def\boV{{\cal V}}
\def\boW{{\cal W}}
\def\cqfd{\hfill$\Box$}
\def\Res{{\,\rm Res}}

\def\Re{{\rm Re}}
\def\Im{{\rm Im}}

\def\disp{\displaystyle}

\newtheorem{theorem}{Theorem}
\newtheorem{proposition}{Proposition}
\newtheorem{definition}{Definition}
\newtheorem{lemma}{Lemma}

\newtheorem{remark}{Remark}

\title{A minimal surface with unbounded curvature}
\author{Martin Traizet}
\begin{document}
\maketitle
\section{Introduction}
\label{section-intro}
The main goal of this paper is to construct a complete, embedded minimal surface
in euclidean space $\R^3$ which has unbounded Gaussian curvature.
Besides being a mathematical curiosity, this example
is also of theoretical interests, as it illustrates some recent
developments in the theory of properly embedded minimal surfaces of infinite
topology.

\medskip

Let me first explain why all examples of complete embedded minimal surfaces known so far
have bounded Gaussian curvature. Almost all known examples fall into one of the
following (non-disjoint) categories :
\begin{enumerate}
\item {\em Finite total curvature minimal surfaces}. Such a surface has a finite number
of catenoidal or planar ends, so it must have bounded curvature.
\item {\em Periodic minimal surfaces which have finite total curvature in the quotient}.
In this case, we know that the quotient is either compact (in the triply-periodic case)
or has a finite number of ends which are asymptotically flat, so again the
curvature is bounded.
This category contains the vast majority of known examples.
\item {\em Properly embedded minimal surfaces with finite genus}.
By a recent result of Meeks, Perez and Ros \cite{mpr4}, such a surface
has bounded curvature. As examples which do not fit in one of the previous
categories, we have the genus one helicoid \cite{howe3} which has one end,
and the Riemann examples with handles constructed in \cite{hauP1} which have
infinitely many planar ends.
\end{enumerate}
As far as I know, the only known examples which do not fit into one of these
categories are the Saddle Towers with infinitely many ends \cite{mrt1} and
the quasi-periodic examples constructed in \cite{matr1}. Both are proven
to have bounded Gaussian curvature (which actually requires some work).

\medskip

In this paper we prove
\begin{theorem}
\label{th1}
There exists a complete, properly embedded minimal surface in euclidean
space $\R^3$ which has unbounded Gauss curvature.
It has infinite genus, infinitely many catenoid type ends, and one limit-end.
\end{theorem}

From the theoretical point of view, the most interesting feature of this
example is its last property.

\medskip

Collin, Kusner, Meeks and Rosenberg \cite{ckmr1} have proven that a properly
embedded minimal surface with infinitely many ends has at most two limit
ends.
Meeks, Perez and Ros \cite{mpr4} have proven that in the finite genus case,
such a surface cannot have one limit end.
The Riemann examples have genus zero and
two limit ends.
However, no example with just one limit end was known, so it seems interesting
to construct an example to illustrate the theory. Of course, it must have infinite genus.

\medskip

Let me point out that the existence of such an example is not completely
unexpected. Indeed, at least heuristically, one can imagine how
to construct one by inductively desingularizing a family of suitable
catenoids. However, we don't have a general enough desingularization theorem
at our disposal yet, and there are fantastic technicalities in trying to carry
out such a construction. So the purpose of this paper is to contruct an example
using another idea, in a somewhat more economical way.

\medskip

Another remark is that if we relax the embeddedness condition, then there are
plenty of known complete, immersed minimal surfaces with unbounded Gaussian
curvature. For example, the example of Nadirashvili \cite{na1} of a
complete minimal immersion in a ball certainly has unbounded curvature.
Embeddedness is a strong constraint on the geometry of minimal surfaces.

\medskip

Heuristically, our example is constructed inductively as follows.
Start with the catenoid and stack a plane on top of it.
Glue a finite number of catenoidal necks in between. After this first
step one gets a Costa Hoffman Meeks surface with three ends. Then iterate this
process infinitely many times, increasing the number of ends by one at each step.
What we need to carry on this construction is a theorem which, from a minimal
surface with $n$ ends, produces a minimal surface with one more end. This
theorem is the main result of this paper and is stated in the next section.

\begin{figure}[h]
\begin{center}
\epsfxsize=10cm
\epsffile{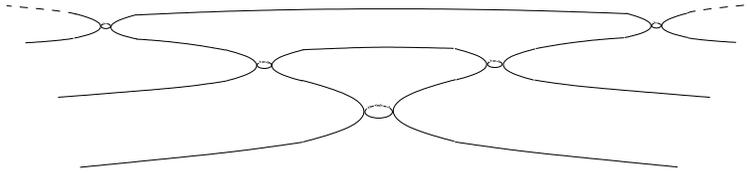}
\caption{A sketch of the surface we get after two steps.
Only two catenoidal necks have been represented at each level for clarity.}
\label{fig1}
\end{center}
\end{figure}

\section{Main result}
\label{section-main}
Given an embedded minimal surface $M$ of finite total curvature in $\R^3$ with
$n$ catenoidal ends, we denote $c_1(M),c_2(M),\cdots,c_n(M)$ the
logarithmic growths of its ends, ordered from bottom to top.
\begin{theorem}
\label{th2}
Let $M$ be an embedded minimal surface in $\R^3$ with finite total curvature,
genus $G$ and $n$ horizontal catenoidal ends with logarithmic growths satisfying
$$c_1(M)<c_2(M)\cdots<c_{n-1}(M)<0<c_n(M).$$
Assume that $M$
has a vertical plane of symmetry and
non-degenerate Weierstrass Representation.

Consider an integer $m\geq 2$ such that
\begin{equation}
\label{eq-condition-sur-m}
m-1>\frac{c_n(M)}{|c_{n-1}(M)|}.
\end{equation}
There exists a one parameter family of embedded minimal surfaces
$(M_t)_{0<t<\varepsilon}$ which has the following
properties:
\begin{enumerate}
\item $M_t$ has finite total curvature, genus $G+m-1$ and $n+1$ catenoidal ends, whose
logarithmic growths satisfy
$$c_1(M_t)<c_2(M_t)\cdots<c_{n}(M_t)<0<c_{n+1}(M_t).$$
\item $M_t$ converges smoothly to $M$ on compact subsets of $\R^3$ when $t\to 0$.
\item $M_t$ has a vertical plane of symmetry and has non-degenerate Weierstrass Representation.
\item The maximum of the absolute value of the Gaussian curvature on $M_t$
is greater than $\disp\frac{(m-1)^2}{2c_n(M)^2}$.
\end{enumerate}
\end{theorem}
This theorem will be proven in 
section \ref{section-construction}.
The definition of ``non-degenerate Weierstrass Representation'' will be given in 
section \ref{section-nondegenerate}.
Roughly speaking, this means that the space of deformations of
$M$, keeping the vertical plane of symmetry, has the expected dimension, namely
$n-1$. In particular, the catenoid has non-degenerate Weierstrass Representation.

\medskip

Heuristically, $M_t$ is constructed by stacking a horizontal plane on top of $M$ and gluing $m$ catenoidal necks placed on a circle in between.
When $t\to 0$, the catenoids drift off to infinity, which is why $M_t$ converges to $M$ on compact subsets of $\R^3$.

\medskip

The catenoidal necks must all have waist radius equal to $r=\disp\frac{c_n(M)}{m-1}$,
they cannot be chosen freely. The Gaussian curvature of such a catenoid along
its waist circle is equal to $-1/r^2$, so this explains the last item of the
theorem (the $1/2$ factor is here only because these catenoids only approximate
the surface).

\medskip

By a simple evaluation of the vertical flux, we obtain, from the size of
the catenoidal necks, the following values for the logarithmic growths of $M_t$
\begin{equation}
\label{eq-log-growths}
\lim_{t\to 0} c_k(M_t)=\left\{
\begin{array}{l}
c_k(M) \mbox{ if $1\leq k\leq n-1$}\\
\frac{-1}{m-1}c_n(M) \mbox{ if $k=n$}\\
\frac{m}{m-1}c_n(M) \mbox{ if $k=n+1$}
\end{array}\right.
\end{equation}
From this we see that the condition \eqref{eq-condition-sur-m} ensures that
$c_{n-1}(M_t)<c_n(M_t)$ as required. The construction is of course possible
without this condition but will not produce an embedded minimal surface.
\section{Proof of theorem 1}
\label{section-limit}
We construct inductively a sequence of minimal surfaces $(S_n)_{n\geq 2}$,
an increasing sequence of balls $(B_n)_{n\geq 2}$ and a sequence of
positive numbers $(C_n)_{n\geq 2}$ with the following
properties:
\begin{enumerate}
\item Each $S_n$ is an embedded
minimal surface of finite total curvature with $n$ catenoidal ends
satisfying
$$c_1(S_n)<c_2(S_n)<\cdots<c_{n-1}(S_n)<0<c_n(S_n)$$
and with a vertical plane of symmetry and non-degenerate Weierstrass
Representation.
\item For all $\ell\geq k\geq 2$, one has
\begin{equation}
\label{eq-*}
k-2<\sup_{S_\ell\cap B_k} |K|<C_k
\qquad \mbox{ and } \qquad
\mbox{Area}(S_\ell\cap B_k)<C_k
\end{equation}
where $K$ denotes the Gaussian curvature.
\end{enumerate}
The process is initiated with $S_2$ equal to the standard catenoid,
$B_2=B(0,2)$ and $C_2$ a suitable constant.
Take $n\geq 2$ and
assume that $S_{\ell}$, $B_{\ell}$ and $C_{\ell}$ have been constructed for
all $\ell\leq n$,
so that \eqref{eq-*} is satisfied for all $2\leq k\leq \ell\leq n$.
We apply theorem \ref{th2} with $M=S_n$ and $m=m_n$ chosen
large enough so that
condition \eqref{eq-condition-sur-m} is satisfied
and $\frac{(m-1)^2}{2 c_n(S_n)^2}>n-1$.
The output of the theorem is a family of minimal surfaces $(M_t)_{0<t<\varepsilon}$
which converges to $S_n$ on each $B_k$ for $k\leq n$. Hence
we can choose $t$ small enough so that
$S_{n+1}=M_t$ satisfies \eqref{eq-*} for all $k\leq n$.
By the last item of theorem \ref{th2}, there are points on $S_{n+1}$ where
$|K|>n-1$. We take a ball $B_{n+1}$ large enough to contain one such point,
and containing $B_n$.
Then we can choose a constant $C_{n+1}$ so that $S_{n+1}$ satisfies \eqref{eq-*}
for $k=n+1$ and we are done.

\medskip

For each $k\geq 2$, the sequence $(S_n\cap B_k)_{n\geq k}$ has uniform
curvature and area estimate, so has a subsequence which converges smoothly
by standard compactness results (theorem 4.2.1 in \cite{pro2}).
By a diagonal process, the sequence $(S_n)_{n\geq 2}$ has a subsequence
which converges smoothly on each $B_k$, to a complete embedded
minimal surface $S_{\infty}$. Now for all $k\geq 2$,
$\sup_{S_{\infty}\cap B_k} |K|\geq k-2$, so $S_{\infty}$ has
unbounded Gaussian curvature and the theorem is proven.
\cqfd
\begin{remark}
All the catenoidal ends of $S_{\infty}$ have negative logarithmic growth.
\end{remark}
In the above argument, we have chosen the sequence $(m_n)_{n\geq 2}$
so that the limit surface $S_{\infty}$ has unbounded Gaussian curvature, but
can we choose it so that $S_{\infty}$ has bounded Gaussian curvature ?

\medskip

The sequence $(m_n)_{n\geq 2}$ must be chosen so that condition \eqref{eq-condition-sur-m}
is satisfied at each step. Using formula \eqref{eq-log-growths}, we have
$$c_{n}(S_{n+1})\simeq \frac{-1}{m_n-1}c_n(S_n)$$
$$c_{n+1}(S_{n+1})\simeq \frac{m_n}{m_n-1}c_n(S_n)$$
where $\simeq$ means that it can be chosen arbitrarily close by taking $t$
small enough.
So condition \eqref{eq-condition-sur-m} reads as
$m_{n+1}-1>m_n$.
Take an arbitrary sequence $(m_n)_{n\geq 2}$ satisfying
\begin{equation}
\label{eq-condition-sur-m_n}
m_2\geq 3
\quad \mbox{ and }\quad
\forall n\geq 2,\quad m_{n+1}\geq  m_n +2
\end{equation}
By the above process, we obtain a sequence of minimal surface $(S_n)_{n\geq 2}$
which converges to an embedded minimal surface $S_{\infty}$
with infinitely many catenoidal ends.
By induction, we have $m_n\geq 2n-1$ and
$$c_n(S_n)\simeq \prod_{i=2}^{n-1}\frac{m_i}{m_i-1}
\leq\prod_{i=2}^{n-1}\left(1+\frac{1}{2i-2}\right)
=O(\sqrt{n}).$$
Hence, $\lim \frac{c_n(S_n)}{m_n-1}=0$.
By the last item of theorem \ref{th2}, this means that whatever the choice
of the sequence $(m_n)_{n\geq 2}$ satisying \eqref{eq-condition-sur-m_n}, the
minimal surface $S_{\infty}$ will have unbounded Gaussian curvature.

Also, we have
$$\forall n\geq 2,\quad c_n(S_{\infty})\simeq \frac{-1}{m_n-1}\prod_{i=2}^{n-1}\frac{m_i}{m_i-1}$$
so depending on the choice of the sequence $(m_n)_{n\geq 2}$, the series
$\sum c_n(S_{\infty})$ can be convergent or divergent.
\section{Non-degenerate Weierstrass Representation}
\label{section-nondegenerate}
Let $M$ be an embedded minimal surface in $\R^3$ with genus $G$
and $n$ horizontal catenoidal ends.
Let $(\Sigma,g,\phi_3)$ be its Weierstrass Representation.
Here $\Sigma$ is a compact Riemann surface, the Gauss map
$g:\Sigma\to\CC=\C\cup\{\infty\}$ is
a meromorphic function and the height differential
$\phi_3$ is a meromorphic 1-form on $\Sigma$ with $n$ simple poles
which we call $q_1,\cdots,q_n$.
These points correspond to the ends of $M$ and are called the punctures.
The degree of the Gauss map is $d=G+n-1$.
Define
$$\phi_1=\frac{1}{2}(g^{-1}-g)\phi_3,\quad
\phi_2=\frac{i}{2}(g^{-1}+g)\phi_3.$$
Our minimal surface $M$ is parametrized on $\Sigma\setminus\{q_1,\cdots,q_n\}$
by
\begin{equation}
\label{eq-weierstrass}
z\mapsto\Re\int_{z_0}^z (\phi_1,\phi_2,\phi_3).
\end{equation}
We assume that $M$ has a vertical plane of symmetry. Without loss of
generality we assume that $M$ is symmetric with respect to the
plane $x_2=0$. On $\Sigma$, this symmetry corresponds to a antiholomorphic
involution $\sigma$ such that $g\circ\sigma=\overline{g}$ and
$\sigma^*\phi_3=\overline{\phi_3}$.
Moreover, $\sigma$ fixes the punctures $q_1,\cdots,q_n$.
\begin{definition}
We say that the triple $(\Sigma,g,\phi_3)$ is $\sigma$-symmetric
if there exists a antiholomorphic involution $\sigma:\Sigma\to\Sigma$
satisfying
$g\circ\sigma=\overline{g}$ and $\sigma^*\phi_3=\overline{\phi_3}$.
\end{definition}

Let us pretend we would like to deform $M$, keeping the vertical plane of
symmetry.
Let us write $(\Sigma_0,g_0,\phi_{3,0})$ for the Weierstrass data
of the minimal surface we are given.
In the following sections, we count how many parameters are
available for $\sigma$-symmetric
deformation of the Weierstrass data
and how many equations need to be solved.
``Non-degenerate Weierstrass Representation''
simply means that the jacobian matrix of equations with
respect to parameters has maximal rank.

\subsection{Deformations of the Weierstrass data}
\label{section-deformation}
Let $\Sigma_0$ be a compact Riemann surface of genus $G$ and
$g_0:\Sigma_0\to\CC$ be a meromorphic function of degree $d$. We see it as a branched
covering of the Riemann sphere, and we would like to parametrize all deformations of $g_0$ by branched coverings $g:\Sigma\to\CC$.
Two branched coverings $g:\Sigma\to\CC$ and $g':\Sigma'\to\CC$ are said to be isomorphic
if there exists a biholomorphic map $\psi:\Sigma\to\Sigma'$ such that
$g=g'\circ\psi$.

In certain cases, the space of isomorphism classes of branched coverings are known to be
smooth complex manifolds. For instance, a degree $d$ branched covering is said to be simple if
each fiber contains at least $d-1$ points (so each fiber contains at most one branch point
and it has branching order 1). The moduli space of simple branched coverings of given degree
is called a Hurwitz space. It is an open complex manifold of dimension $2G+2d-2$. The list of branching
values of $g$ provide local coordinates on this space. More generally, the moduli space of coverings of
degree $d$ and with $n$ branch points
form a smooth complex manifold of dimension $n$
\cite{Fulton}.

Now if our covering has a branch point of branching order $k\geq 2$, when deforming it,
this branch point may split into several smaller order branch points,
whose branching orders sum up to $k$.
It is not true anymore that the list of branching values provide local coordinates
(see the example at the beginning of appendix \ref{appendixB}).
The moduli space of branched covering of given degree, with no restriction on the branch points, is not a smooth complex manifold.
It is singular at those coverings who have a branch point which is fixed by
a nontrivial automorphism of the covering.

In the case of minimal surfaces, the punctures $q_1,\cdots,q_n$ are
$n$ distinguished points on $\Sigma$ at which $g$ takes the value alternately
$0$ and $\infty$ with multiplicity one.
We call the data $(\Sigma,g,q_1,\cdots,q_n)$ a
{\em marked covering}. We say that two marked coverings
$(\Sigma,g,q_1,\cdots,q_n)$ and
$(\Sigma,g',q'_1,\cdots,q'_n)$ are isomorphic if there exists a
biholomorphic map $\psi:\Sigma\to\Sigma'$ such that
$g=g'\circ\psi$ and $\psi(x_i)=x'_i$ for all $1\leq i\leq n$.

The moduli space of marked coverings has
a structure of complex manifold of dimension $2G+2d-2$.
If, moreover, we require the coverings to be $\sigma$-symmetric,
then it is a real manifold of dimension $2G+2d-2$.
We present in appendix \ref{appendixB} one way to define local coordinates on this space.

As a conclusion, we can parametrize all $\sigma$-symmetric deformations 
of the marked covering
$(\Sigma_0,g_0,q_1,\cdots,q_n)$ as
$(\Sigma_a,g_a,q_1(a),\cdots,q_n(a))$ with a parameter $a\in \R^{2G+2d-2}$
in a neighborhood of $0$.
On each $\Sigma_a$ we have a antiholomorphic involution
$\sigma$ such that $g_a\circ\sigma=\overline{g_a}$.
The function $g$ takes alternately the value $0$ and $\infty$ at $q_1(a),\dots,q_n(a)$,
with multiplicity one. These points are fixed by $\sigma$.

Then we would like to write all candidates for the height differential $\phi_3$ on
$\Sigma_a$. It needs simple poles at the punctures $q_1,\cdots,q_n$.
The most natural way to define such a meromorphic 1-form is to prescribe
residues and periods.
Consider a canonical homology basis $A_1,\cdots,A_G,B_1,\cdots,B_G$ of $\Sigma_a$
(depending continuously on $a$).
Given complex parameters $c=(c_1,\cdots,c_{n-1})\in\C^{n-1}$ and $\alpha=(\alpha_1,\cdots,\alpha_G)\in\C^G$
we define $\phi_3$ as the unique meromorphic 1-form on $\Sigma_a$
with simple poles at $q_1,\cdots,q_n$ and the following
residues and $A$-periods :
$$\Res_{q_i}\phi_3=-c_i \qquad 1\leq i\leq n-1,$$
$$\int_{A_i}\phi_3=2\pi i \alpha_i \qquad 1\leq i\leq G.$$
The residue at $q_n$ is $c_1+\cdots+c_{n-1}$ by the Residue theorem.
We would like $\phi_3$ to be $\sigma$-symmetric.
This translates into simple conditions on the parameters $c$ and $\alpha$,
provided we choose carefuly the homology basis.
\begin{definition}
Let $\Sigma$ be a Riemann surface with a antiholomorphic involution
$\sigma$.
We say that a canonical homology basis
$\{A_1,\cdots,A_G,B_1,\cdots,B_G\}$
is $\sigma$-symmetric if there exists an action of
$\sigma$ on $\{1,\cdots,G\}$ such that for all $1\leq i\leq G$,
$\sigma(A_i)\sim -A_{\sigma(i)}$ and $\sigma(B_i)\sim B_{\sigma(i)}$.
\end{definition}
One can construct a $\sigma$-symmetric canonical homology basis
of $\Sigma=\Sigma_a$ as follows.
The fixed set of
$\sigma$ divides $\Sigma$ into two components $\Sigma^+$
and $\Sigma^-$, so that $\sigma(\Sigma^+)=\Sigma^-$.
Let $G'$ be the genus of $\Sigma^+$
and $k$ the number of its boundary components. Then $G=2G'+k-1$.
We consider $2G'$ cycles $A_1,\cdots,A_{G'},B_1,\cdots,B_{G'}$ in $\Sigma^+$
such that the only non-zero intersection numbers are $A_i.B_i=1$.
For $1\leq i\leq G'$, we define $A_{G'+i}=-\sigma(A_i)$ and
$B_{G'+i}=\sigma(B_i)$.
We take $B_{2G'+1},\cdots,B_{2G'+k-1}$ to be the
homology classes of all boundary components of $\Sigma^+$ but one.
For $1\leq i\leq k-1$, we choose a cycle $A_{2G'+i}$ which intersects
$B_{2G'+i}$ (with intersection number $1$) and the remaining component
of the boundary of $\Sigma^+$.

The condition $\sigma^*\phi_3=\overline{\phi_3}$ is equivalent to
$\alpha_{\sigma(i)}=\overline{\alpha_i}$, for $1\leq i\leq G$, and
$c_i\in\R$ for $1\leq i\leq n-1$.
The parameters $c_i$ correspond geometrically to the logarithmic growths
of the ends.
The condition on $\alpha$ defines a space of real dimension $G$.
\begin{remark} Of course the Period Problem will require the parameters
$\alpha_i$ to be real numbers, but it is better to leave this as an equation
to be solved, else the dimension of the parameter space depends on $G'$.
\end{remark}
At this point we have defined a family of triples $(\Sigma_a,g_a,\phi_3)$ depending
on the parameter $X=(a,c,\alpha)$. The parameter space
has real dimension $5G+3n-5$.
We write $X^0$ for the value of $X$ which gives the Weierstrass data
$(\Sigma_0,g_0,\phi_{3,0})$ of the minimal
surface $M$ we were given.
In general we will not write explicitely the dependance of objects on parameters,
so we will write $(\Sigma,g,\phi_3)$, beeing understood that everything depends on the
parameters.
\subsection{The equations}
\label{section-equations}
In order to be define an immersed minimal surface, the triple $(\Sigma,g,\phi_3)$
must satisfy the following conditions :
\begin{enumerate}
\item At any zero of the height differential $\phi_3$, the Gauss map $g$ needs a zero or
a pole, with the same multiplicity. We call this the zero/pole equation.
\item At each puncture $q_1,\cdots,q_n$, the residues of $\phi_1$ and $\phi_2$ must
be real.
\item For all $\nu=1,2,3$ and $1\leq i\leq G$ we need
$\disp\Re\int_{A_i}\phi_{\nu}=\Re\int_{B_i}\phi_{\nu}=0$.
\end{enumerate}
Point 2 and 3 guarantee that \eqref{eq-weierstrass} is well defined, and
point 1 that it is an immersion.
Let us count how many equations we have to solve, taking into account the $\sigma$-symmetry.
\begin{itemize}
\item Regarding point 2, the residues of $\phi_1$ are real and the residues of $\phi_2$ are
imaginary by symmetry. Also, the residue of $\phi_2$
at $q_j$ is $\frac{i}{2}\Res\, g^{-1}\phi_3$ if $g(q_j)=0$, and
$\frac{i}{2}\Res\, g\phi_3$ if $g(q_j)=\infty$.
Provided point 1 is satisfied, the only poles of $g\phi_3$ and $g^{-1}\phi_3$ are at the punctures.
Applying the Residue Theorem to $g\phi_3$ and $g^{-1}\phi_3$,
it suffices to solve the equation $\Res_{q_j}\phi_2=0$ for $1\leq i\leq n-2$.
Therefore, point 2 counts as $n-2$ real equations.
\item Regarding the Period Problem, we have the symmetries $\sigma^*\phi_{\nu}=(-1)^{\nu+1}
\overline{\phi_{\nu}}$, $\nu=1,2,3$. From this we obtain, for $1\leq i\leq G$,
$$\Re\int_{A_{\sigma(i)}}\phi_{\nu}=(-1)^{\nu}\Re\int_{A_i}\phi_{\nu}$$
$$\Re\int_{B_{\sigma(i)}}\phi_{\nu}=(-1)^{\nu+1}\Re\int_{B_i}\phi_{\nu}.$$
So the period problem reduces to $3G$ real equations. Explicitely, for the canonical basis
that we defined in the previous section, these equations are
$$\Re\int_{A_i}\phi_{\nu}=0,\qquad 1\leq i\leq G',\quad \nu=1,3,$$
$$\Re\int_{B_i}\phi_{\nu}=0,\qquad G'+1\leq i\leq G,\quad \nu=1,3,$$
$$\Re\int_{A_i}\phi_2=0,\qquad G'+1\leq i\leq G,$$
$$\Re\int_{B_i}\phi_2=0,\qquad 1\leq i\leq G'.$$
\item Regarding point 1, the height differential $\phi_3$ has $2G+n-2$ zeros, counting
multiplicity. Let us first assume that the zeros of $\phi_{3,0}$ are simple. Then
this remains true when the parameter $X$ is close to $X^0$, and we may label
them $\zeta_1(X),\cdots,\zeta_{2G+n-2}(X)$ so that they depend continuously on $X$.
If $g_0$ has a zero (resp. a pole) at $\zeta_i(X^0)$, we need to solve the equation $g_a(\zeta_i(X))=0$
(resp. $g(a)^{-1}(\zeta_i(X))=0$). Taking into account the $\sigma$-symmetry (which implies that the
set of zeros of $\phi_3$ is invariant by $\sigma$), these are $2G+n-2$ real equations (the equations
are either real or conjugate by pairs).
\item In case $\phi_{3,0}$ has a zero of multiplicity $k\geq 2$ at some point $\zeta$, the problem is
that this zero may split into several zeros of smaller multiplicity when we deform the Weierstrass
data. Let us assume for example that $g_0$ has a zero (of multiplicity $k$) at $\zeta$.
By the very definition of $(\Sigma_a,g_a)$, there is a local coordinate $z$ on $\Sigma_a$ such that
$$g_a(z)=z^k+\sum_{i=0}^{k-1} a_i z^i$$
where the coefficients $a_i$ are components of the parameter vector $a$ (and actually $a_{k-1}=0$).
By the Weierstrass Preparation Theorem, we can write
$$\phi_3=f(X,z)(z^k+\sum_{i=0}^{k-1} b_i(X)z^i)dz$$
where the coefficients $b_i(X)$ are holomorphic functions of $X$.
The polynomial $z^k+\sum b_i z^i$ is called the {\em Weierstrass Polynomial}
of $\phi_3$, it depends on the choice of the coordinate $z$.
We solve the equations $b_i(X)=a_i$, $0\leq i\leq k-1$.
Taking into account the $\sigma$-symmetry we have the same number of equations as before.
\end{itemize}
Everything together, we have to solve $5G+2n-4$ real equations.
Let us write these equations as $\boF(X)=0$.
Recall that the dimension of the parameter space is $5G+3n-5$.
\begin{definition}
\label{def-non-degenerate}
We say that $M$ has non-degenerate Weierstrass Representation if the
differential of $\boF$ at $X^0$ is onto. Equivalently, its kernel
has dimension $n-1$.
\end{definition}
Example : the catenoid has non-degenerate Weierstrass Representation,
because in this case, $G=0$ and $n=2$, so
the parameter space has dimension 1 and there are no equations to solve.
As expected, the catenoid comes into a one-parameter family (by scaling).

As a straightforward application of the implicit function theorem, we have
\begin{proposition}
If $M$ has non-degenerate Weierstrass
Representation, then it can be deformed into a family of non-congruent
minimal surfaces depending on $n-1$ parameters, all symmetric
with respect to the plane $x_2=0$.
\end{proposition}
\begin{remark}\em
We do not require that the parameters in this deformation are the logarithmic
growths of the ends, so this notion of non-degeneracy should be weaker than
the usual non-degeneracy defined in term of Jacobi fields \cite{pro4}, altough
I have not been able to prove that one of them implies the other.
\end{remark}
\section{Proof of theorem 2}
\label{section-construction}
\def\wtSigma{\widetilde{\Sigma}}
\def\wtg{\widetilde{g}}
\def\wtG{\widetilde{G}}
We construct the family of minimal surfaces $(M_t)_{0<t<\varepsilon}$
by writing down candidates
for its Weierstrass Representation and
then solving the zero/pole equation and the Period Problem.
We define the Riemann surface and the Gauss map by opening nodes.
The height differential $\phi_3$
is defined by prescribing periods and residues.
The equations are solved using the implicit function theorem, using
the fact that $M$ has non-degenerate Weierstrass 
Representation.
\subsection{Opening nodes}
\label{section-opening-nodes}
Since $M$ has non-degenerate Weierstrass Representation, we are given
in particular
a family of branched coverings $(\Sigma,g)=(\Sigma_a,g_a)$
depending on the parameter $a$.
Let $q_n=q_n(a)\in\Sigma_a$ be the point which corresponds to the top end.
Without loss of generality we may assume that $g(q_n)=0$.

We consider two copies of the complex plane,
denoted $\C^-$ and $\C^+$.
We choose $m$ distinct, non-zero points $p_1^-,\cdots,p_m^-$ in $\C^-$ and
$m$ distinct points $p_1^+,\cdots,p_m^+$ in $\C^+$.
Identify the point $q_n$ with the point $0$ in $\C^-$.
For $1\leq i\leq m$, identify the point $p_i^-$ with the point $p_i^+$.
This defines a Riemann surfaces with $m+1$ nodes which we call
$\wtSigma_0$.

We define two meromorphic functions $g^-$ on $\C^-$ and $g^+$ on $\C^+$ by
$$g^-(z)=\frac{\beta_0^-}{z}+\sum_{i=1}^m \frac{\beta_i^-}{z-p_i^-}$$
$$g^+(z)=\sum_{i=1}^m\frac{\beta_i^+}{z-p_i^+}.$$
Here $\beta_0^-,\cdots,\beta_m^-$ and $\beta_1^+,\cdots,\beta_m^+$
are non-zero complex parameters.
We write $\beta^{-}=(\beta^{-}_0,\cdots,\beta^{-}_m)$,
$\beta^{+}=(\beta^{+}_1,\cdots,\beta^{+}_m)$
and $p^{\pm}=(p_1^{\pm},\cdots,p_m^{\pm})$.

Since $g$ has a simple zero at $q_n$ and $g^{\pm}$ have simple poles, we
can fix a small number $0<\epsilon<1$ such that $v_0^-:=g$ is a diffeomorphism
from a small neighborhood $U_0^-$ of $q_n$ to the disk $D(0,\epsilon)$,
$v_0^+:=1/g^-$
is a diffeomorphism from a small neighborhood $U_0^+$ of $0$ in $\C^-$
to $D(0,\epsilon)$ and for each $i=1,\cdots,m$, $v_i^{\pm}:=1/g^{\pm}$ is
a diffeomorphism from a small neighborhood $U_i^{\pm}$ of $p_i^{\pm}$ in $\C^{\pm}$
to $D(0,\epsilon)$.
We use $v_i^{\pm}$ as local complex coordinates to open nodes.

Consider a real parameter $t$ such that $0<t<\epsilon^2$.
We remove the disk $|v_0^-|<\frac{t}{\epsilon}$ from $U_0^-$ and
the disk $|v_0^+|<\frac{t}{\epsilon}$ from $U_0^+$.
We identify the point $z\in U_0^-$ with the point $z'\in U_0^+$ such that
$v_0^-(z) v_0^+(z')=t$.
This is equivalent to 
$$g(z)=t g^-(z').$$
For each $1\leq i\leq m$, we remove the disks
$|v_i^-|<\frac{t^2}{\epsilon}$ from $U_i^-$ and
$|v_i^+|<\frac{t^2}{\epsilon}$ from $U_i^+$.
We identify the point $z\in U_i^-$ with
the point $z'\in U_i^+$ such that
$v_i^-(z)v_i^+(z)=t^2$.
This is equivalent to
$$ tg^-(z)=\frac{1}{t g^+(z)}.$$
This defines a Riemann surface. We compactify
it by adding the points at infinity in $\C^-$ and
$\C^+$ and call it $\wtSigma$.
Its genus is $\wtG=G+m-1$.

By a slight abuse of language, we will denote by $\Sigma\subset\wtSigma$,
$\C^-\subset\wtSigma$ and $\C^+\subset\wtSigma$ the domains
$\Sigma$, $\C^-$ and $\C^+$ minus the disks that were removed when
opening nodes.

We define the Gauss map $\wtg$ on $\wtSigma$ by
$$\wtg(z)=\left\{\begin{array}{ll}
g(z) & \mbox{ if $z\in\Sigma$}\\
tg^-(z) & \mbox{ if $z\in \C^-$}\\
\disp\frac{1}{tg^{+}(z)} & \mbox{ if $z\in\C^+$}
\end{array}\right.$$
This is a well defined meromorphic function on $\wtSigma$ because
$\wtg(z)=\wtg(z')$ whenever $z$ and $z'$ are identified.

Next we would like $(\wtSigma,\wtg)$ to have the required symmetry.
By hypothesis, $\Sigma$ does have a antiholomorphic involution
$\sigma$ such that $g\circ\sigma=\overline{g}$.
Let $V$ be the linear subspace of $\C^{m-1}$ defined by
$z_{m-i}=\overline{z_i}$ for $1\leq i\leq m-1$.
This is a real vector space of dimension $m-1$.
We require that $\beta^-\in\R\times V\times\R$, and
$\beta^+$, $p^+$ and $p^-$ all
belong to $V\times\R$.
Under these assumptions, we have
$g^-(\overline{z})=\overline{g^-(z)}$ and $g^+(\overline{z})=\overline{g^+(z)}$.
We define $\widetilde{\sigma}:\wtSigma\to\wtSigma$ by
$\widetilde{\sigma}=\sigma$ on $\Sigma$ and
$\widetilde{\sigma}(z)=\overline{z}$ on $\C^-$ and $\C^+$.
Then thanks to the fact that $t$ is real, $\widetilde{\sigma}$
is well defined on $\wtSigma$, and $\wtg\circ\widetilde{\sigma}
=\overline{\wtg}$.
We will write $\sigma$ for $\widetilde{\sigma}$ as no confusion can occur.
\subsection{The height differential}
\label{section-height-differential}
\def\wtphi{\widetilde{\phi}}
As in section \ref{section-deformation}, we define the height differential on $\wtSigma$
by prescribing periods and residues, so we need to define a canonical
homology basis of $\wtSigma$.
The cycles $A_1,\cdots,A_G,B_1,\cdots,B_G$ on $\Sigma$ define us $2G$ cycles on $\wtSigma$.
For $1\leq i\leq m-1$,
let $A_{G+i}$ be the homology class of the circle
$C(p_i^+,\epsilon)$ with the positive orientation. This circle is homologous
to the circle $C(p_i^-,\epsilon)$ with the negative orientation.

For $1\leq i\leq \frac{m}{2}$, we define $B_{G+i}$ as the composition of the
following paths :
\begin{enumerate}
\item a path from the point $v_m^+=-\epsilon$ to the point $v_i^+=\epsilon$ in
$\C^+$,
\item the segment from $v_i^+=\epsilon$ to $v_i^+=\frac{t^2}{\epsilon}$,
\item a path from the point $v_i^-=\epsilon$ to the point
$v_m^-=-\epsilon$ in $\C^-$,
\item the segment from $v_m^-=-\epsilon$ to $v_m^-=-\frac{t^2}{\epsilon}$.
\end{enumerate}
In point 1 and 3, the curves must depend continuously on parameters, and must avoid
all disks around the nodes. In particular, if $i=\frac{m}{2}$, we cannot take
the straight segment in point 3 because it goes through the origin (which has
been removed from $\C^-$).
Also, this definition only makes sense for $t\neq 0$.
For $\frac{m}{2}<i\leq m-1$, we define $B_{G+i}$ as $\sigma(B_{G+m-i})$.
Then $A_1,\cdots,A_{\wtG},B_1,\cdots,B_{\wtG}$ is a $\sigma$-symmetric
canonical homology basis of $\wtSigma$.

Let $\infty^-$ and $\infty^+$ denote the point at infinity in $\C^-$ and $\C^+$.
The punctures (corresponding to the $n+1$ catenoidal ends) are at $q_1,\cdots,q_{n-1},\infty^-$ and
$\infty^+$.
We define the height differential $\wtphi_3$ on $\wtSigma$ as in section \ref{section-deformation} by prescribing its
$A$-periods and its residues at all punctures but one. Actually, by the residue theorem, 
prescribing the residue at $\infty^+$ is the same as prescribing the period on the circle
$C(p_m^+,\epsilon)$. So we define $\wtphi_3$ on $\wtSigma$ as the unique meromorphic
1-form with simple poles at the punctures with the following residues and periods : 
$$\int_{A_j}\wtphi_3=2\pi i \alpha_j, \qquad 1\leq j\leq G$$
$$\Res_{q_i}\wtphi_3=-c_i,\qquad 1\leq i\leq n-1,$$
$$\int_{C(p_j^+,\epsilon)}\wtphi_3=2\pi i \gamma_j, \qquad 1\leq j\leq m$$
The parameters $c=(c_1,\cdots,c_{n-1})$ and $\alpha=(\alpha_1,\cdots,\alpha_{G})$
are as in section \ref{section-deformation}.
The parameter $\gamma=(\gamma_1,\cdots,\gamma_m)$ is in the vector space $V\times\R$
defined above,
so that $\wtphi_3$ satisfies $\sigma^*\wtphi_3=\overline{\wtphi_3}$.
By the residue theorem in $\C^+$ and $\wtSigma$, we have
\begin{equation}
\label{eq-residue+}
\Res_{\infty^+}\wtphi_3=-\sum_{i=1}^m\gamma_j,
\end{equation}
\begin{equation}
\label{eq-residue-}
\Res_{\infty^-}\wtphi_3=\sum_{j=1}^{n-1}c_j+\sum_{j=1}^m\gamma_j.
\end{equation}
By standard results \cite{fay,masur}, on 
a Riemann surface with nodes,
the notion of holomorphic (or meromorphic)
1-form must be replaced by that of a regular differential, which means that it has simple
poles on each side of each node, with opposite residues. In the case at hand, this means that
when $t=0$, $\wtphi_3$ is meromorphic in $\Sigma$, $\C^-$ and $\C^+$, with simples poles
at $q_1,\cdots,q_n$ in $\Sigma$, simple poles at $0$, $p_1^-,\cdots,p_m^-$ and $\infty^-$
in $\C^-$ and simples poles at $p_1^+,\cdots,p_m^+$ and $\infty^+$ in $\C^+$.
The above period and residue prescription define $\wtphi_3$ uniquely,
and moreover $\wtphi_3$
depends analytically (away from the poles) on $t$ and all other parameters.
Given the residues and periods of $\wtphi_3$, we have, when $t=0$ :
$$\wtphi_3=\left\{\begin{array}{l}
\phi_3 \qquad \mbox{ in $\Sigma$},\\
\disp-\sum_{i=1}^{n-1} c_i\frac{dz}{z}-\sum_{i=1}^m\frac{\gamma_i}{z-p_i^-}dz
\qquad\mbox{ in $\C^-$},\\
\disp\sum_{i=1}^m\frac{\gamma_i}{z-p_i^+}dz\qquad\mbox{ in $\C^+$}.
\end{array}\right.$$
The key point is that when $t=0$, the restriction of $\wtg$ and $\wtphi_3$ to $\Sigma$ are
$g$ and $\phi_3$ as defined in section \ref{section-deformation}.
This allow us to use the non-degeneracy hypothesis.
\subsection{Central value of the parameters}
\label{section-central-value}
All the parameters that we have introduced vary in a neighborhood of a central value, which is
the point at which we will apply the implicit function theorem.

The central value of the $t$ parameter is zero. The central value of the parameters
$a$, $\alpha_1,\cdots,\alpha_G$ and $c_1,\cdots,c_{n-1}$ are the values such that
$(\Sigma,g,\phi_3)=(\Sigma_0,g_0,\phi_{3,0})$ is the Weierstrass Representation of our given minimal surface $M$.

The central value of the other parameters will be found by solving the equations, but
it may help the reader to give it here.
Without loss of generality, we may assume by scaling that the logarithmic growth of the top end of $M$ is $c_n=1$.
The central value of the parameters $\beta^+$, $\beta^-$ and $\gamma$
are then given by
$$\beta_0^-=-1,\qquad
\beta_i^+=\beta_i^-=\gamma_i=\frac{1}{m-1}, \qquad 1\leq i\leq m.$$
The central value of the parameters $p^{\pm}$ are given by
$$p_i^+=\omega^{-i},p_i^-=\omega^{i},\qquad 1\leq i\leq m$$
where $\omega=e^{2\pi i/m}$ is a primitive $m$-th root of unity. 

When all parameters have their central value, we have
$$g^-(z) =\frac{-1}{z}+\frac{1}{m-1}\sum_{i=1}^m \frac{1}{z-\omega^i}
=\frac{-1}{z}+\frac{mz^{m-1}}{(m-1)(z^m-1)}$$
$$g^+(z) =\frac{1}{m-1}\sum_{i=1}^m \frac{1}{z-\omega^i}
=\frac{mz^{m-1}}{(m-1)(z^m-1)}.$$
Moreover, $\wtphi_3=-g^-(z) dz$ in $\C^-$ and $\wtphi_3=g^+(z)dz$
in $\C^+$.
\subsection{The equations}
\label{section-equations2}
As in section \ref{section-equations},
we have to solve the following equations
to ensure that the Weierstrass data $(\wtSigma,\wtg,\wtphi_3)$ defines
a minimal immersion.
\begin{itemize}
\item[1a.] At each zero of $\wtphi_3$ in $\Sigma\subset\wtSigma$, the gauss map
$\wtg=g$ needs a zero or a pole, with the same multiplicity.
\item[1b.] At each zero of $\wtphi_3$ in $\C^{\pm}$, $g^{\pm}$ needs a zero with
the same multiplicity.
\item [2a.] At each puncture $q_1,\cdots q_{n-2}$, the residues of
$\phi_1$ and $\phi_2$ must be real.
\item [2b.] The residues of $\phi_1$ and $\phi_2$ at $\infty^+$ must be real.
\item [3a.] $\disp \Re\int_{A_i}\wtphi_{\nu}=\Re\int_{B_i}\wtphi_{\nu}=0$
\quad \mbox{ for }$1\leq i\leq G$, $\nu=1,2,3$.
\item [3b.] $\disp \Re\int_{A_{G+i}}\wtphi_{\nu}=\Re\int_{B_{G+i}}\wtphi_{\nu}=0$
\quad \mbox{ for }$1\leq i\leq m-1$, $\nu=1,2,3$.
\end{itemize}
In the following points, we study how each of these equations extend to $t=0$.
\begin{itemize}
\item Let $\boF$ be the collection of the equations in point 1a, 2a and 3a.
When $t=0$, the restriction of $\wtg$ and $\wtphi_3$ to $\Sigma$ are
$g$ and $\phi_3$, so these equations are exactly the same as the ones in
section \ref{section-equations}. So when $t=0$, the function $\boF$ equals the function
$\boF$ defined in section \ref{section-equations}.
The non-degeneracy hypothesis will take
care of this equation.
\item Regarding point 1b, at the central value of the parameters, $\wtphi_3$
has $m$ simple zeros in $\C^-$, so this remains true for nearby values of
the parameters. We may call these zeros $\zeta_1,\cdots,\zeta_m$ so that
$\sigma(\zeta_i)=\zeta_{m+1-i}$.
Let $\boZ^-=(g^-(\zeta_1),\cdots,g^-(\zeta_m))$
(the letter $\boZ$ stands for ``zero'').
This is an analytic function of all parameters.
Moreover, $\boZ^-_{m+1-i}=\overline{\boZ^-_i}$ so $\boZ^-$ takes value in a real space
of dimension $m$.
\item At the central value, $\wtphi_3$ has one zero of multiplicity $m-1$ at the origin.
We may write $g^+(z)=\frac{P(z)}{Q(z)}$ where $P$ is a unitary polynomial of degree
$m-1$ whose coefficients are (polynomial) functions of the parameters $\beta_i^+$ and
$p_i^+$. Let $R$ be the Weierstrass polynomial of $\wtphi_3$ in a neighborhood of $0$.
We define $\boZ^+=P-R$, this is a real polynomial of degree $m-2$ whose coefficients are
analytic functions of all parameters.
We need to solve the equation $\boZ^+=0$.
\item Regarding the $A$-periods in point 3b, we define
$$\boV^A_j=\Re\int_{C(p_j^+,\epsilon)}\wtphi_3=-2\pi \Im(\gamma_j),\qquad 1\leq j\leq m-1,$$
$$\boH^A_j=\frac{1}{t}\left(\Re\int_{C(p_j^+,\epsilon)}\wtphi_1+i\Re\int_{C(p_j^+,\epsilon)}\wtphi_2\right),
\qquad 1\leq j\leq m.$$
(The letters $\boV$ stand for ``vertical'', $\boH$ for ``horizontal'' and $A$ for
``$A$-cycles''.)
Note that the equation $\boH_A=0$ takes care of point 2b by the Residue Theorem.
The symmetry gives $\boH^A_{m-i}=-\overline{\boH^A_i}$, so $\boH^A$ takes value in a space
of real dimension $m$.
The function $\boH^A$ extends analytically to $t=0$ by the following computation :
\begin{eqnarray}
\boH^A_i &=& \frac{1}{2t}\left(\overline{\int_{C(p_i^+,\epsilon)}\wtg^{-1}\wtphi_3}
-\int_{C(p_i^+,\epsilon)}\wtg\wtphi_3\right)
\nonumber\\
&=& \frac{1}{2t}\left(\overline{\int_{C(p_i^+,\epsilon)}\wtg^{-1}\wtphi_3}+\int_{C(p_i^-,\epsilon)}
\wtg\wtphi_3\right)
\nonumber\\
&=&\frac{1}{2}\left(\overline{\int_{C(p_i^+,\epsilon)}g^+\wtphi_3}
+\int_{C(p_i^-,\epsilon)} g^-\wtphi_3\right)
\label{eq-Horiz-A}
\end{eqnarray}
\item By lemma 1 in \cite{tra1}, the function $\int_{B_{G+i}}\wtphi_3-(\gamma_i-\gamma_m)\log t^2$
extends to an analytic function of all parameters at $t=0$. We make the change of variable
$t=\exp(\frac{-1}{\tau^2})$
where $\tau$ is a real parameter in a neighborhood of zero. We define
the renormalised vertical $B$-periods as
$$\boV^B_i=\tau^2\Re\int_{B_{G+i}}\wtphi_3,\qquad 1\leq i\leq m-1.$$
Then $\boV^B=(\boV^B_1,\cdots,\boV^B_{m-1})$ extends as a smooth function of all parameters
at $\tau=0$, with value
$$\boV^B_i=-2\Re(\gamma_i-\gamma_m)\qquad \mbox{ at $\tau=0$.}$$
The symmetry gives $\boV^B_{m-i}=\boV^B_i$.
\item We define the renormalised horizontal $B$-periods as
$$\boH^B_j=t\left(\Re\int_{B_{G+j}}\wtphi_1+i\Re\int_{B_{G+j}}\wtphi_2\right),\qquad
1\leq j\leq m-1$$
By lemma 2 in \cite{tra1}, $\boH^B=(\boH^B_1,\cdots,\boH^B_{m-1})$ extends to a smooth function
of all parameters at $\tau=0$ with value
\begin{equation}
\label{eq-Horiz-B}
\boH^B_i=\frac{1}{2}\overline{\int_{p_i^-}^{p_m^-}\frac{\wtphi_3}{g^-}}
-\frac{1}{2}\int_{p_m^+}^{p_i^+}\frac{\wtphi_3}{g^+}\mbox{ at $\tau=0$}.
\end{equation}
The symmetry gives $\boH^B_{m-i}=\overline{\boH^B_i}$ so $\boH^B$ takes value
in the space $V$.
\end{itemize}
\subsection{Solving the equations}
\label{section-solving}
\def\dotgamma{\dot{\gamma}}
\def\dotbeta{\dot{\beta}}
\def\dotp{\dot{p}}
\def\wtboF{\widetilde{\cal F}}
Let $X$ be the collection of all parameters but $\tau$.
We denote by $X^0$ the central value of the parameters.
Let $\wtboF(\tau,X)$
be the collection of the equations that we have to solve, namely
$\wtboF=(\boF,\boZ^-,\boZ^+,\boV^A,\boH^A,\boV^B,\boH^B)$.
We want to solve the equation $\wtboF(\tau,X)$ to get $X$ as an implicit function
of $\tau$.
\begin{lemma}
\label{lemma1}
We have $\wtboF(0,X^0)=0$, and the partial differential of $\wtboF$
with respect to $X$ at $(0,X^0)$ is onto.
\end{lemma}
Proof : we make a change of parameters so that the partial differential is triangular by
blocks.
Let
$$\gamma_i=\gamma_m+\dotgamma_i, \qquad 1\leq i\leq m-1,$$
$$\beta_i^-=\gamma_i+\dotbeta_i^-,\qquad 1\leq i\leq m,$$
$$\beta_0^-=\sum_{i=1}^{n-1}c_i+\dotbeta_0^-,$$
$$\beta_i^+=\gamma_i+\dotbeta_i^+,\qquad 1\leq i\leq m,$$
$$p_i^+=\overline{p_i^-}+\dotp_i^+,\qquad 1\leq i\leq m.$$
We write
$\dotgamma=(\dotgamma_1,\cdots,\dotgamma_{m-1})\in V$,
$\dotbeta^-=(\dotbeta^-_0,\cdots,\dotbeta^-_m)\in\R\times V\times\R$,
$\dotbeta^+=(\dotbeta^+_1,\cdots,\dotbeta^+_m)\in V\times\R$
and $\dotp^+=(\dotp^+_1,\cdots,\dotp^+_m)\in V\times\R$.
The central value of each of these new parameters is $0$.
Now the parameters are
$a$,
$\alpha$,
$c$,
$\gamma_m$,
$\dotgamma$,
$\dotbeta^-$,
$\dotbeta^+$,
$p^-$,
$\dotp^+$ and $\tau$.
In the following points,
we evaluate the partial differential at the central value
of each equation with respect to all these parameters except $\tau$.
At the same time we check that the equations are satisfied at the central value
of the parameters.
\begin{itemize}
\item By the non-degeneracy hypothesis, {\em the partial differential of $\boF$ with respect to the
parameters $(a,\alpha,c)$ is onto. Moreover, its partial derivative with respect to all other
parameters is zero.}
\item {\em  The partial differential of $\boZ^-$ with respect to $\dotbeta^-$ is onto,
and all other partial derivatives of $\boZ^-$ are zero.}

Indeed, when $t=0$, the zeros of $\wtphi_3$ do not depend on $\beta^-$ anymore, so $\boZ^-$ is a linear function
of $\dotbeta^-$. If $\dotbeta^-$ is in the kernel of the partial differential, then $\boZ^-(\dotbeta^-)=0$,
so $\wtphi_3$ and $g^- dz$ have the same zeros in $\C^-$. Since they have the same poles, they are proportionnal,
$\wtphi_3=\lambda g^- dz$. Hence the kernel has dimension 1. Since the $\dotbeta^-$ space has dimension $m+1$
and the target space has dimension $m$, the partial differential is onto.
The second statement holds because if $\dotbeta^-=0$, then $\wtphi_3=-g^-dz$ in $\C^-$, so they have the same zeros
hence $\boZ^-=0$.

\item In the exact same way, {\em the partial differential of $\boZ^+$ with respect to $\dotbeta^+$ is onto
and all other partial derivatives of $\boZ^+$ are zero.}

\item {\em The partial differential of $(\boV^A,\boV^B)$ with respect to $\dotgamma$ is an isomorphism. All other
partial derivatives are zero.}

Indeed, in term of the new parameters we have
$\boV^A_i=-2\pi\Im(\dotgamma_i)$
and $\boV^B_i=-2\Re(\dotgamma_i)$, so the partial derivative is injective. Because of
the symmetry, the domain and target spaces have the same dimension, namely $m-1$.

\item {\em The partial differential of $\boH^B$ with respect to $\dotp^+$ is onto. The only other nonzero partial
derivatives are the partial differentials with respect to $\dotbeta^+$ and $\dotbeta^-$.}   

Indeed, if $\dotbeta^+=0$ and $\dotbeta^-=0$, then $\wtphi_3=-g^- dz$ in $\C^-$ and $\wtphi_3=g^+dz$ in
$\C^+$. By equation \eqref{eq-Horiz-B}, we get $\boH^B_i=\frac{1}{2}(\dotp_m^+-\dotp_i^+)$.
The statement readily follows.

\item {\em The partial differential of $\boH^A$ with respect to $(p_1^-,\cdots,p_{m-1}^-,\gamma_m)$ is an isomorphism.}

To prove this, assume that all parameters but $p^-$ and $\gamma_m$ have their central value. 
Then $\wtphi_3=-g^- dz$ in $\C^-$ and $\wtphi_3=g^+ dz$ in $\C^+$, so
formula \eqref{eq-Horiz-A} gives
\begin{eqnarray*}
\boH^A_i&=&-\pi\overline{\Res_{p_i^+}(g^+)^2}
-\pi\Res_{p_i^-}(g^-)^2\\
&=&-2\pi\sum_{j\neq i}\frac{\gamma_m^2}{\overline{p_i^+}
-\overline{p_j^+}}-2\pi\sum_{j\neq i}\frac{\gamma_m^2}{
p_i^- -p_j^-}+2\pi\frac{\gamma_m}{p_i^-}\\
&=&-4\pi\sum_{j\neq i}\frac{\gamma_m^2}{p_i^- -p_j^-}
+2\pi\frac{\gamma_m}{p_i^-}.
\end{eqnarray*}
This implies that
$$\sum_{i=1}^m p_i^-\boH^A_i
=-2\pi m(m-1)\gamma_m^2+2\pi m\gamma_m.$$
When $\gamma_m$ has its central value, namely $\frac{1}{m-1}$, the right hand side is zero. When
$p_1^-,\cdots,p_m^-$ have their central value, all terms in the left sum are equal by symmetry,
so all are zero. This proves that $\boH^A=0$ at the central value.

Now consider the matrix of the partial differential of $\boH^A$ with respect to
$(p_1^-,\cdots,p_{m-1}^-,\gamma_m)$. Perform the row operation
$R_m\to \sum p_i^- R_i$. By
the previous formula, we obtain a matrix of the form
$$\disp\left(\begin{array}{cc} A & \cdot \\ 0 & -2\pi m\end{array}\right)$$
where $A$ is a square matrix of order $m-1$.
This matrix is proven to be invertible in appendix \ref{appendixA}, which proves the
statement.
\end{itemize}
The lemma readily follows from these statements (the matrix of the partial differential
has block triangular form).
\cqfd
\subsection{Proof of Theorem 2}
\label{section-embedded}
\def\wtpsi{\widetilde{\psi}}
By lemma \ref{lemma1} and the implicit function theorem, for $\tau$ in a neighborhood of $0$, there
exists a smooth function $X(\tau)$ such that $\wtboF(\tau,X(\tau))=0$.
For $t>0$ close to zero,
let use write $(\wtSigma_t,\wtg_t,\wtphi_{3,t})$ for the Weierstrass data
corresponding to the value $\tau=(\log t)^{-1/2}$ and $X=X(\tau)$ of the
parameters.
We choose a base point $z_0\in\Sigma$, then this Weierstrass data defines us a
minimal immersion $\wtpsi_t$ on $\wtSigma_t$ minus the punctures. 
Let $M_t$ be its image.
In the following points we prove that the family $(M_t)_{0<t<\varepsilon}$
has all the properties claimed in Theorem \ref{th2}.
\begin{itemize}
\item $M_t$ has $n+1$ catenoidal ends at $q_1,\cdots,q_{n-1}$,
$\infty^-$ and $\infty^+$. The logarithmic growths are the opposite of the
residue of $\wtphi_3$ at these points, so by equations \eqref{eq-residue+}
and \eqref{eq-residue-},
their limit value when $t\to 0$
are
$c_1,\cdots,c_{n-1},1-\frac{m}{m-1}$ and $\frac{m}{m-1}$.
Since we have scaled $M$ so that $c_n(M)=1$, this gives formula \eqref{eq-log-growths}.
\item {\em $M_t$ converges to $M$ on compact subsets of $\R^3$.}

This follows from the fact that $\wtg_t$ converges to $g_0$
on $\Sigma$ and $\wtphi_{3,t}$
converges to $\phi_{3,0}$ on compact subsets of $\Sigma$ minus the punctures,
(We translate $M$ so that the image of $z_0$ is the origin.)

\item {\em $M_t$ has non-degenerate Weierstrass Representation.}

Indeed, since having
maximal rank is an open property, the differential of $\wtboF$ at
$(\tau,X(\tau))$ remains onto for $\tau$ close to $0$.
The only issue here is that the parameters $t$, $p^{\pm}$,
$\beta^{\pm}$ for $(\wtSigma,\wtg)$ are not the
right ones for the definition of non-degeneracy. However, the coordinates
that we use on the Hurwitz space are analytic functions of the parameters
$t,p^{\pm},\beta^{\pm}$, so if we denote by $\widehat{X}$ the parameters that we
use in the definition of non-degeneracy,
we have $\widehat{X}=f(\tau,X)$ for
some smooth map $f$. Write the equations that we have to solve for non-degeneracy
as $\widehat{\boF}(\widehat{X})=0$. Then $\wtboF=\widehat{\boF}\circ f$,
so the differential of $\widehat{\boF}$ is onto.
\item {\em $M_t$ is embedded.}

To prove this statement, we study the asymptotic behavior of $\wtpsi_t$ on
each of the domains $\Sigma$, $\C^-$ and $\C^+$ when $t\to 0$.
On $\C^+$ we have 
$$\lim_{t\to 0} t\wtphi_1=-\frac{dz}{2},
\quad\lim_{t\to 0} t\wtphi_2=i\frac{dz}{2},
\quad\lim_{t\to 0}\wtphi_3
=\frac{1}{m-1}\sum_{i=1}^m\frac{dz}{z-\omega^i}.$$
Define $\widehat{\psi}_t$ on $\C^+$
as the composition of $\wtpsi_t-\wtpsi_t(0)$ with the affine transformation
$(x_1,x_2,x_3)\mapsto (-2t x_1,-2t x_2,x_3)$. Then
$$\lim_{t\to 0}\widehat{\psi}_t(z)=(\Re\, z,\Im\, z, u^+(z))$$
where $u^+$ is the harmonic function
$$u^+(z)=\frac{1}{m-1}\sum_{i=1}^m \log|z-\omega^i|.$$
So the image of $\widehat{\psi}_t$ converges to the graph of $u^+$.
For $h$ large enough, the graph of $u^+$ intersects the plane $x_3=-h$
in $m$ closed convex curves, so the same is true for the image
of $\widehat{\psi}_t$ for $t$ small enough.
As a conclusion, we can find a height $c_1$ (depending on $t$) such that the image
of $\C^+$ by $\wtpsi_t$
cuts the plane $x_3=c_1$ in $m$ closed convex curves.
We call $S^+$ the part which is above this plane. The surface $S^+$ is
embedded (as a graph) and has one upward-going catenoidal end.

In the same way, after horizontal
scaling by $-2t$ and vertical translation, the image of
$\C^-$ by $\widetilde{\psi}$ converges to the graph of $u^-(\overline{z})$,
where
$$u^-(z)=-\log |z| + \frac{1}{m-1}\sum_{i=1}^m \log|z-\omega^i|.$$
For $h$ large enough, the graph of $u^-$ intersects the plane $x_3=h$ in
$m$ closed convex curves and the plane $x_3=-h$ in two closed convex curves,
one inside the other. 
Again we may find some heights $c_2$ and $c_3$, with $c_3<c_2<c_1$ such that 
for $t$ small enough, the image of $\C^-$ by $\wtpsi_t$ cuts the plane
$x_3=c_2$ in $m$ closed convex curves and the plane $x_3=c_3$ in two closed
convex curves, one inside the other. Let $S^-$ be the part bounded by the $m$
top curves and the inside bottom curve. It is an embedded surface with one
downward catenoidal end.

Finally, since the top end of $M$ is catenoidal, we may find some height $c_4<c_3$ such
that the image of $\Sigma$ by $\wtpsi_t$ cuts the plane $x_3=c_4$ in one closed
convex curves (and what is above is an annulus). Let $S$ be the part which is
below this plane. It is embedded because $M$ is. 

The pieces $S$, $S^-$ and $S^+$ are disjoint.
(For $S^-$ and $S$, this uses the maximum principle and the fact that the logarithmic
growth of the end of $S^-$ is larger than the logarithmic growth of the top end
of $S$).
Each component of the complementary set in $M_t$ of $S\cup S^-\cup S^+$
is a minimal annulus bounded by two closed
convex curves in parallel planes. By a theorem of Shiffman \cite{sh1},
such an annulus is fibered by horizontal curves. It follows that
$M_t$ is embedded.
\begin{figure}[h]
\begin{center}
\epsfxsize=7cm
\epsffile{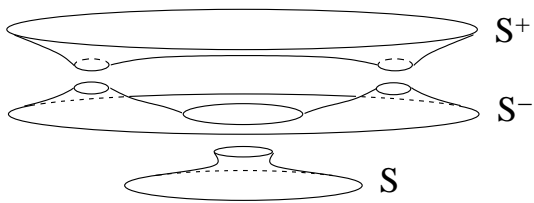}
\caption{The pieces $S^+$, $S^-$ and $S$ (in case $M$ is a catenoid).}
\label{fig2}
\end{center}
\end{figure}
\item {\em For $t$ small enough,
the maximum of $|K|$ on $M_t$ is greater than $\frac{1}{2}(m-1)^2$.}

Indeed, assume by contradiction that this is not true. Then we can find a sequence
$(t_n)_n$ converging to zero such that the Gaussian curvature on $M_{t_n}$ is
bounded by $\frac{1}{2}(m-1)^2$. Let $A_n$ be the image of the annulus
bounded by the circles $|v_1^+|=\epsilon$ and $|v_i^-|=\epsilon$ in $M_{t_n}$.
Translate $A_n$ so that the point where the Gauss map is one is at the origin.
The images of the boundary circles are close to circles of radius $O(\frac{1}{t_n})$,
so for any ball $B(0,R)$, $A_n\cap B(0,R)$ is properly embedded in $B(0,R)$
for $n$ large enough.
As we have uniform Gaussian curvature and area estimate (by the monotonicity
formula) for $A_n$, a subsequence of $(A_n)_n$
converges smoothly on compact subsets of $\R^3$ to a complete embedded minimal annulus
(\cite{pro2}, theorem 4.2.1) 
hence a catenoid. As the flux of $A_n$ converges to $(0,0,\frac{2\pi}{m-1})$, the
limit catenoid has waist radius $\frac{1}{m-1}$. Since the maximum of the Gaussian
curvature on this catenoid is $(m-1)^2$, we have a contradiction.
\end{itemize}
\appendix
\section{A matrix computation}
\label{appendixA}
\begin{lemma}
\label{lemma-CHM}
Consider an integer $m\geq 2$.
Let $\omega=e^{2\pi i/m}$. Define $p_i=\omega^i$ for $1\leq i\leq m$.
Consider the order $m-1$ square matrix $A$
defined by
$$a_{ii}=\frac{m-1}{(p_i)^2}-\sum_{1\leq j\leq m\atop j\neq i}
\frac{2}{(p_i-p_j)^2}$$
$$a_{ij}=\frac{2}{(p_i-p_j)^2}\quad \mbox{ if $j\neq i$}$$
Then $A$ is invertible.
\end{lemma}
Proof. The proof relies on the following observation :
if $z\in\C$ is such that $|z|=1$, then
$$1-2\Re\frac{1}{(1+z)^2}=\frac{2}{|1+z|^2}$$
To see this, observe that the transformation $z\to\frac{1}{z+1}$ maps
the unit circle to the line $\Re(z)=\frac{1}{2}$. If $w$ is on this line,
an elementary computation gives $1-2 \Re (w^2)=2 |w|^2$.

We prove that $A$ has dominant diagonal.
We have
$$(p_i)^2 a_{ii}=\sum_{j=1}^{m-1}\left(1-\frac{2}{(1-\omega^j)^2}\right).$$
Since this is a real number,
$$(p_i)^2 a_{ii}=\sum_{j=1}^{m-1}\Re\left(1-\frac{2}{(1-\omega^j)^2}\right)
=\sum_{j=1}^{m-1}\frac{2}{|1-\omega^j|^2}$$
$$|a_{ii}|=\sum_{j=1}^{m-1}\frac{2}{|1-\omega^j|^2}
=\sum_{1\leq j\leq m-1\atop j\neq i} |a_{ij}|
+\frac{2}{|p_i-1|^2}.$$
Hence $|a_{ii}|>\disp\sum_{j\neq i}|a_{ij}|$ so $A$ is invertible
by Hadamard theorem.
\section{Deformations of a branched covering}
\label{appendixB}
Let us start with an example to illustrate the problem.
We consider the following two deformations of the covering $z\mapsto z^4$
of the Riemann sphere :
$$f_t(z)=z^4+4tz^3,$$
$$g_t(z)=z^4+4t\alpha z^3+4t^2\alpha^2 z^2,
\qquad \mbox{ with $\alpha^4=-27$}.$$
By explicit computations, they have the same branching
values close to $0$, namely $0$, $0$ and $-3t$.
They are not isomorphic because
$f_t$ has a branch point of branching order 2 at the origin,
whereas $g_t$ has three simple branch points
(at $0$, $-t\alpha$ and $-2t\alpha$).
The conclusion of this example is that in general one cannot parametrize
the deformations of a covering with a high order branch point by their
branching values.

\begin{definition}
Fix a value $y_0$ in the Riemann sphere. 
A {\em marked branched covering} is a triple $(\Sigma,g,x)$ where $\Sigma$
is a compact Riemann surface, $g$ is a meromorphic function on $\Sigma$
and $x\in\Sigma$ is a regular point of $g$ (i.e. not a branch point) such
that $g(x)=y_0$.

Two marked branched coverings $(\Sigma,g,x)$ and $(\Sigma',g',x')$ are
isomorphic if there exists a biholomorphic map $\psi:\Sigma\to\Sigma'$ such that
$g=g'\circ\psi$ and $\psi(x)=x'$.
\end{definition}
Here we assume for simplicity that there is just one marked point, but there could
be several, this does not change anything.

Let $(\Sigma_0,g_0,x_0)$ be a marked branched covering, $G$ the genus of $\Sigma_0$
and $d$ the degree of $g_0$.
We first construct an explicit family of deformations depending on $2G+2d-2$ complex parameters.

Let $p_1,\cdots,p_r$ be the branch points of $g_0$. For $1\leq i\leq r$, let $k_i\geq 2$ be the
multiplicity of $g_0$ at $p_i$ (so the branching order is $k_i-1$) and $q_i=g(p_i)$ the branching
value. By using a Moebius transformation, we may assume without loss of generality that
all branching values are finite.

We may find $\varepsilon>0$ small enough so that the following holds: for each 
$1\leq i\leq r$, $g_0$ is a branched covering of degree $k_i$ from a closed topological disk $U_i$ containing $p_i$
to the closed disk $\overline{D}(q_i,\varepsilon)$. Moreover, the disks $U_1,\cdots,U_r$ are disjoint, and two
disks $\overline{D}(q_i,\varepsilon)$ and
$\overline{D}(q_j,\varepsilon)$ are either disjoint or equal.
We choose, on the boundary of each disk $U_i$, a point $x_i$ such that $g_0(x_i)=q_i+\varepsilon$, which we will use as a marking.

For each $1\leq i\leq r$, consider the polynomial
$$h_i=z^{k_i}+q_i+\sum_{j=0}^{k_i-2} a_{ij} z^i.$$
If the complex parameters $a_{ij}$ are small enough, all branching values (except $\infty$) of $h_i$ are
inside the disk $D(q_i,\varepsilon)$, and
$h_i$ is a branched covering of degree $k_i$ from a closed disk
$V_i$ to the closed disk $\overline{D}(q_i,\varepsilon)$. We mark this covering with the point
$y_i\in\partial V_i$ such that $h_i(y_i)=q_i+\varepsilon$ which is
closest to $\varepsilon^{1/k_i}$. It depends continuously on the parameters
$a_{ij}$ provided they are small enough.

Then we remove the interior of
$U_i$ from $\Sigma_0$ and replace it by $V_i$, indentifying a point
$z\in\partial U_i$ with the point $z'\in\partial V_i$ such that $h_i(z')=g_0(z)$.
Of course, there
are $k_i$ possible choices for $z'$,
but the markings allow us to solve this indetermination. More
precisely, we may parametrise the boundary of $U_i$ as $c(t)$, $t\in[0,2k_i\pi]$, so that
$c(0)=x_i$ and $g_0(c(t))=q_i+\varepsilon e^{it}$. We may parametrise the boundary of $V_i$ as
$c'(t)$, $t\in[0,2k_i\pi]$, so that $c'(0)=y_i$ and $h_i(c'(t))=q_i+\varepsilon e^{it}$.
We identify the point $c(t)$ with the point $c'(t)$ for all $t$.

Doing this for all $1\leq i\leq r$ defines a compact Riemann surface which we denote
$\Sigma_{a}$, where $a$ denotes the collection of all parameters $a_{ij}$, $1\leq i\leq r$,
$0\leq j\leq k_i-2$.
Observe that the number of parameters is precisely the total branching order
of $g_0$, namely $2G+2d-2$.
On $\Sigma_{a}$ we may define a meromorphic function $g_a$ by $g_a=g_0$ on $\Sigma_0$ minus the
disks $U_i$, and $g_a=h_i$ on each disk $V_i$.
Finally, we mark the covering $(\Sigma_a,g_a)$ with the point $x_0$ we were given (which may be
assumed to be outside the disks $U_i$).

\begin{remark}\em If $k_i=2$, then $h_i(z)=z^2+q_i+a_{i0}$ and the branching value of
$g_a$ in $V_i$ is $q_i+a_{i0}$. So in the case of simple branch points, our parameters
are (up to translation by $q_i$) the branching values, as in the standard
parametrisation of the Hurwitz spaces.
\end{remark}
\begin{remark}\em The symmetric functions of the
branching values of $g_a$ inside 
each $V_i$ are polynomial functions of the numbers $a_{ij}$,
because the symmetric functions of the branching values of a polynomial
may be expressed in function of its coefficients.
\end{remark}

\begin{remark}\em
The family $(\Sigma_a)_a$ is a holomorphic family of compact Riemann surfaces, in
the sense that there exists a complex manifold $\boW$ of dimension 
$2G+2d-1$ and a holomorphic map $\pi:\boW\to \C^{2G+2d-2}$ such that
for $a$ close to $0$, $\pi^{-1}(a)$ is isomorphic to $\Sigma_a$.
\end{remark}

\begin{remark}\em
If $(\Sigma_0,g_0,x_0)$ is $\sigma$-symmetric,
one can define an action of $\sigma$
on $\{1,\cdots,r\}$ by $\sigma(p_i)=p_{\sigma(i)}$. Choose the marking
$x_i$ so that $x_{\sigma(i)}=\sigma(x_i)$.
Then $(\Sigma_a,g_a,x_0)$
is $\sigma$-symmetric if and only if $a$ satisfies
$a_{\sigma(i)j}=\overline{a_{ij}}$ for all $1\leq i\leq r$ and $0\leq j\leq k_i-2$.
Indeed, if $a$ satisfies this condition, then one can define an involution
$\widetilde{\sigma}$ on $\Sigma_a$ by $\widetilde{\sigma}=\sigma$ on
$\Sigma_0$ minus the disks $U_i$, and for each $1\leq i\leq r$,
$\widetilde{\sigma}$ maps the point $z\in V_i$ to the point $\overline{z}$
in $V_{\sigma(i)}$. It is straightforward to check that $\widetilde{\sigma}$
is well defined and satisfies $g_a\circ\widetilde{\sigma}=\overline{g_a}$.
The ``only if'' part of the statement is a consequence of uniqueness,
proposition \ref{prop2} below.
\end{remark}
\medskip

We want to prove that with this construction, we obtain all deformations
of $(\Sigma_0,g_0)$, up to isomorphism of marked coverings.
We need the following
\begin{lemma}
Let $U\subset\C$ be a closed disk and $g:U\to \overline{D}(0,1)$ be a branched holomorphic covering of degree $k$.
The exists a closed disk $V\subset\C$, a biholomorphic mapping $\psi : U\to V$ and a
polynomial $h$ of degree $k$
such that $V=h^{-1}(\overline{D}(0,1))$ and $g=h\circ\psi$.
Moreover, $h$ is unique up to composition (on the right) by a transformation
of the form $z\mapsto \alpha z+\beta$.
\end{lemma} 
Proof : choose a point $x_0$ in $\partial U$ such that $g(x_0)=1$. We may parametrize the boundary
of $U$ as $c(t)$, $t\in[0,2\pi]$ so that $g(c(t))=e^{kit}$ and $c(0)=x_0$.
We glue together $U$ and $\CC\setminus D(0,1)$
by identifying the point $c(t)$ with the point $e^{it}$. We obtain a genus zero compact Riemann surface
$\widetilde{U}$ on which we can define a meromorphic function $\widetilde{g}$ by $\widetilde{g}=g$ on
$U$ and $\widetilde{g}=z^k$ on $\CC\setminus D(0,1)$.
By the uniformization theorem, there exists a biholomorphic mapping $\psi:\widetilde{U}\to \CC$.
We may choose $\psi$ so that $\psi(\infty)=\infty$.
Let $h=g\circ \psi^{-1}$. Then $h$ is a meromorphic function on $\CC$ with a single pole of
multiplicity $k$ at infinity, so it is a polynomial of degree $k$.
We restrict $\psi$ to $U$
and let $V=\psi(U)$, the first statement of the proposition is proven.

Regarding uniqueness, assume that we have two degree $k$
polynomials $h$ and $h'$ such that $D(0,1)$ contains all the branching
values of $h$ and $h'$ (except infinity) and the branched coverings
$h:V\to \overline{D}(0,1)$ and $h':V'\to\overline{D}(0,1)$ are isomorphic
by $\psi:V\to V'$. Then $h:\C\setminus V\to \C\setminus \overline{D}$
and $h':\C\setminus V'\to\C\setminus\overline{D}$ are isomorphic
(unbranched) coverings, so we may extend $\psi$ into $\widetilde{\psi}:
\C\to\C$ such that $h=h'\circ\widetilde{\psi}$.
Then $\widetilde{\psi}$ is an automorphism
of the Riemann sphere which maps $\infty$ to $\infty$,
so $\widetilde{\psi}(z)=\alpha z+\beta$.
\cqfd

Next we prove that different choices of the parameters $a$ give non-isomorphic
marked coverings.
\begin{proposition}
\label{prop2}
If the branched coverings $(\Sigma_a,g_a,x_0)$ and $(\Sigma_{a'},g_{a'},x_0)$
are isomorphic then $a=a'$.
\end{proposition}
For example, let us consider the branched covering $z\mapsto z^3$ and the following
two deformations : $z\mapsto z^3+tz$ and $z\mapsto z^3+tjz$, where $j$ is a primitive
cubic root of unity. As branched coverings, they are isomorphic by
$z\mapsto j^2 z$. As marked coverings (with the marking chosen as above)
they are not isomorphic.

Proof of the proposition :
let $\psi:\Sigma_a\mapsto\Sigma_{a'}$ be the covering isomorphism.
Let $\Omega$ be the Riemann sphere minus the disks $D(q_i,\varepsilon)$.
Then by construction, $g_a^{-1}(\Omega)=g_{a'}^{-1}(\Omega)=g_0^{-1}(\Omega)$.
So the restriction of $\psi$ to $g_0^{-1}(\Omega)$ is an automorphism of
the (unbranched) covering $g_0:g_0^{-1}(\Omega)\to\Omega$. Since
$\psi(x_0)=x_0$, it must be the identity. In particular, $\psi(x_i)=x_i$
for all $1\leq i\leq r$, so $\psi$ preserves all the markings.

Then for any $1\leq i\leq r$, consider the polynomials $h_i:V_i\to D(q_i,\varepsilon)$
and $h'_i:V'_i\to D(q_i,\varepsilon)$ which are used to construct respectively
$\Sigma_a$ and $\Sigma_{a'}$. Then $h_i=h'_i\circ\psi$ on $V_i$, so by the 
uniqueness part of the lemma, the restriction of $\psi$ to $V_i$ has the form
$z\mapsto \alpha z+\beta$. As $h_i$ and $h'_i$ are unitary polynomials, we must
have $\alpha^{k_i}=1$. Since they have no term of degree $k_i-1$, we must
have $\beta=0$. Since $\psi$ preserves the markings, $\alpha=1$ by our choice
of the markings of the polynomials. Hence $h_i=h'_i$ for all $1\leq i\leq r$,
which means that $a=a'$.
\cqfd

By a deformation of $(\Sigma_0,g_0)$, we mean the following. Let $S$ be $\Sigma_0$ seen
as a differentiable 2-manifold, i.e. forgetting the conformal structure. Let
$(g_t)_{0\leq t\leq 1}$ be a continuous family of branched coverings from $S$ to $\CC$.
It is well known that each $g_t$ induces a conformal structure on $S$ which we denote
by $\Sigma_t$. We also assume that each $g_t$ is close enough to $g_0$ in the sense that
all the branching values of $g_t$ remain inside the union of the disks $D(q_i,\varepsilon)$.
Then we say that each $(\Sigma_t,g_t)$ is a deformation of $(\Sigma_0,g_0)$.
If $(\Sigma_0,g_0)$ is marked with the point $x_0$, then each $(\Sigma_t,g_t)$ admits a unique
marking $x_0(t)$, depending continuously on $t$, such that $g_t(x_0(t))=g_0(x_0)$.
\begin{proposition}
If $(\Sigma,g)$ is a deformation of $(\Sigma_0,g_0)$, there exists
$a\in\C^{2G+2d-2}$ such that $(\Sigma,g)$ is isomorphic to $(\Sigma_a,g_a)$.
\end{proposition}
Proof : let $\Omega$ be the complement in the Riemann sphere of the disks $D(q_i,\varepsilon)$.
Then the restrictions
$g_0:g_0^{-1}(\Omega)\to\Omega$ and $g:g^{-1}(\Omega)\to\Omega$ are isomorphic (unbranched) coverings  
so there exists an isomorphism $\psi:g_0^{-1}(\Omega)\to g^{-1}(\Omega)$ such that
$g_0=g\circ\psi$.  
We can easily extend $\psi$ to those components of $\Sigma_0\setminus g_0^{-1}(\Omega)$ in which
$g_0$ has no branch point. Indeed, if $U$ is such a component, let us call $U'$ the component
of $\Sigma\setminus g^{-1}(\Omega)$ which is bounded by $\psi(\partial U)$. Then $g_0$ and
$g$ are diffeomorphisms
from respectively $U$ and $U'$ to the same disk $D(q_i,\varepsilon)$,
so we may extend $\psi$ as $(g|_{U'})^{-1}\circ g_0$ on $U$.

Finally, let us consider one of the disks $U_i$ in which $g_0$ has a branch point of mutiplicity $k_i$.
Let $U'_i$ be the component of $\Sigma\setminus g^{-1}(\Omega)$ which is bounded by
$\psi(\partial U_i)$. Then $U'_i$ is a disk and $g:U'_i\to D(q_i,\varepsilon)$ is a branched covering
of degree $k_i$ (altough the branch point may have split into several ones).
Recall that we have a marking $x_i\in\partial U_i$, and let us call
$x'_i=\psi(x_i)\in \partial U'_i$.
By the lemma, there exists a polynomial $h_i$ of degree $k_i$
and a diffeomorphism $\psi_i : U'_i\to V_i$ such that $g=h_i\circ \psi_i$ on $U'_i$.
By composing $h_i$ on the right by a transformation $z\mapsto \alpha z+\beta$,
we may assume that $h_i$ is a unitary polynomial with no term of degree
$k_i-1$, and moreover, $\psi_i(x'_i)$ is the solution of $h_i(z)=q_i+\varepsilon$
which is closest to $\varepsilon^{1/k_i}$.
Let $(\Sigma_a,h_a)$ be the branched covering constructed with the polynomials
$h_1,\cdots,h_r$ we have found.
Define $\widetilde{\psi}:\Sigma_a\to \Sigma$ by 
$\widetilde{\psi}=\psi$ on $\Sigma_0$ minus the disks $U_1,\cdots,U_r$
and $\widetilde{\psi}=\psi_i^{-1}$ on each $V_i$. Then $\widetilde{\psi}$
is an isomorphism between the marked coverings $(\Sigma_a,g_a)$
and $(\Sigma,g)$.
\cqfd

\end{document}